\documentclass[12pt]{amsart}
\usepackage{amsmath,amsthm,amssymb}

\newcommand{\al}{\alpha}

\newcommand{\fr}{\mathcal{F}}

\newcommand{\ity}{\infty}
\newcommand{\C}{\mathbb{C}}

\newcommand{\N}{\mathbb{N}}

\numberwithin{equation}{section}
\newtheorem{theorem}{Theorem}[section]
\newtheorem{lemma}[theorem]{Lemma}
\newtheorem{corollary}[theorem]{Corollary}
\newtheorem*{thmA}{Theorem A}
\newtheorem*{thmB}{Theorem B}
\newtheorem*{FFT}{First Fundamental Theorem}
\newtheorem*{LDT}{Logarithmic Derivative Lemma}
\theoremstyle{remark}
\newtheorem{remark}[theorem]{Remark}
\newtheorem{example}[theorem]{Example}

\makeatletter
\@namedef{subjclassname@2010}{%
  \textup{2010} Mathematics Subject Classification}
\makeatother

\begin{document}

\title[Normality Criteria for Families of Meromorphic Functions ]{Normality Criteria for Families of Meromorphic Functions}

\thanks{The research work of the second author is supported by research fellowship from CSIR India.}
\author[S. Kumar	]{Sanjay Kumar}
\address{Department of Mathematics, Deen Dayal Upadhyaya College, University of Delhi,
Delhi--110 078, India }
\email{sanjpant@gmail.com}
\author[P. Rani]{Poonam Rani}
\address{Department of Mathematics, University of Delhi,
  Delhi--110 007, India} \email{pnmrani753@gmail.com}

\begin{abstract}
In this paper we prove some normality criteria for a family of meromorphic functions, which involves the zeros of certain differential polynomials generated by the members of the family.
\end{abstract}



 \maketitle

\section{Introduction and main results}
 An important aspect  of the theory of complex analytic functions is to find normality criteria for families of meromorphic functions. The notion of normal families was introduced by Paul Montel in 1907. Let us begin by recalling the definition. A family   of meromorphic (holomorphic) functions defined on a domain $D\subset \C$ is said  to be normal in the domain, if every sequence in the family  has a subsequence which converges spherically  uniformly on compact subsets of  $D$ to a meromorphic (holomorphic) function or to $\ity$.\\

 According to Bloch's principle, a family of meromorphic functions in a domain $D$ possessing  a property which reduces a meromorphic function in the plane to a constant, makes the family normal in the domain $D$. Although Bloch's principle is not true in general, many authors established normality criteria for families of meromorphic functions having such properties. \\

 Hayman \cite{Hay1} proved a result which states that: {\it{Let $n\geq 5$ be a positive integer and $a(\neq 0), b$ be two finite complex numbers. If a meromorphic function $f$ in $\C$ satisfies, $f'+af^n \neq b,$  then $f$ is  a constant.}} The normality related to this result was conjectured by Hayman. Confirming this conjecture of Hayman, Drasin ~\cite{DD} proved a normality criterion  for holomorphic function which says: {\it{Let $\fr$ be a family of holomorphic functions in the unit disc $\Delta$, and for a fixed $n\geq3$ and $a(\neq0), b$ are finite complex numbers suppose that for each $f\in \fr$, $f'+af^n\neq b$ in $\Delta$. Then $\fr$ is normal. }} Schwick ~\cite{Sch1} also proved a normality criterion which confirms the following result: {\it {Let $k, n$ be positive integers with $n\geq k+4$ and $a, b\in \C$ with $a\neq0$. Let $\fr$ be a family of meromorphic functions in a domain $D$. If for each $f\in \fr$, $f^{(k)}-af^n\neq b$, then $\fr $ is normal.}}\\

 Similar to above results Fang and Zalcman ~\cite{FZ} proved a normality criterion which states: {\it Let $\fr $ be a family of meromorphic functions in a domain $D$, let $n\geq2$ be a positive integer, and let $a(\neq 0), b \in \C$. If for each  $f\in\fr,$ all zeros of $f$ are multiple and $f+a(f')^n\neq b$ on $D$, then $\fr$ is normal on $D$.} Yan Xu et al ~\cite{XWL} extended this result and proved: {\it Let $\fr $ be a family of meromorphic functions in a domain $D$, let $k, n(\geq k+1)$ be two positive integers, and let $a(\neq 0), b \in \C$. If for each  $f\in\fr,$ all zeros of $f$ are of multiplicity at least $k+1$ and $f+a(f^{(k)})^n\neq b$ on $D$, then $\fr$ is normal on $D$.} They conjectured that the conclusion of their result still holds for $n\geq2.$ Confirming their conjecture, C. L. Lei et al ~\cite{LFZ} proved {\it Let $\fr $ be a family of meromorphic functions in a domain $D$, let $n\geq2$ be a positive integer, and let $a(\neq 0), b \in \C$. If for each  $f\in\fr,$ all zeros of $f$ are of multiplicity at least $k+1$ and $f+a(f^{k})^n\neq b$ on $D$, then $\fr$ is normal on $D$.}\\

 On the other side, Xia and Xu ~\cite{XX} proved the following result:\\

 \begin{thmA}Let $\fr$ be a family of meromorphic functions defined on a domain $D\subset \C.$ Let $\psi(\not\equiv0)$ be a holomorphic function in $D$, and $k$ be a positive integer. Suppose that for every function $f\in\fr$, $f\neq 0$, $f^{(k)}\neq 0$ and all zeros of $f^{(k)}-\psi(z)$ have multiplicities at least ${\frac{k+2}{k}}$. If for $k = 1,$ $\psi$ has only zeros with multiplicities at most 2, and for $k \geq 2,$ $\psi$  has only simple zeros, then F is normal in D.
 \end{thmA}

 Recently, B. Deng et al ~\cite{D} studied the case where $f^{(k)}-af^n-b$ has some distinct zeros in $D$ and they proved the following results:\\

   \begin{thmB}Let $\fr $ be a family of meromorphic functions in a domain $D$, let $m, n, k$ be three positive integers such that $n\geq k+m+2$, and let $a(\neq 0), b \in \C$. If for each  $f\in\fr,$ all zeros of $f$ are of multiplicity at least $k$ and $f^{(k)}-af^n- b$ has at most $m$ distinct zeros, then $\fr$ is normal on $D$.
\end{thmB}

 It is natural to ask  if  $f^{(k)}-af^n-\psi$  has zeros with some multiplicities for a holomorphic function $\psi(z)\neq0$. In this paper we investigate this situation and prove the following results:
 \begin{theorem}\label{thm1}
Let $\fr$ be a family of meromorphic functions in a domain $D$. Let $d, k \ \text{and\ } m$  be positive integers such that  ${\frac{k+1}{d}+ \frac{k+2}{2m}<\frac{1}{2}}$. Let   $a(z)$ and $\psi(z)(\neq0)$ be  holomorphic functions in $D$. If for each $f\in\fr$, all zeros of $f$ are of multiplicity at least $d$, poles of $f$ are multiple and for a positive integer $n$, all zeros of $f^{(k)}(z)-a(z)f^n(z)-\psi(z)$ are of multiplicity at least $m$ in $D$, then $\fr$ is normal in $D$.
\end{theorem}
For a family of holomorphic functions we have the following strengthened version:
\begin{theorem}\label{thm2}
Let $\fr$ be a family of holomorphic functions in a domain $D$. Let $d, k,\ \text{and\ } m$ be  positive integers such that ${\frac{k+1}{d}+ \frac{1}{m}<1}$.  Let $a(z)$ and $\psi(z)(\neq0)$ be  holomorphic functions in $D$. If for each $f\in\fr$, all zeros of $f$ are of multiplicity at least $d$  and for a positive integer $n$, all zeros of $f^{(k)}(z)-a(z)f^n(z)-\psi(z)$ are of multiplicity at least $m$ in $D$, then $\fr$ is normal in $D$.
\end{theorem}
\begin{corollary}
Let $\fr$ be a family of holomorphic functions on a domain $D$, let $k, m, n$ be three positive integers. If for each $f\in \fr$, $f$ has zeros of multiplicities at least $m$, all zeros of $f^{(k)}-1$  are of multiplicity at least $n$ in $D$ and ${\frac{k+1}{m}+ \frac{1}{n}<1}$, then $\fr$ is normal in $D$.
\end{corollary}
The following example shows that $n$ can't be zero in Theorem \ref{thm2}.
\begin{example}
Let $D:=\{z\in \C:|z|<1\}$. Consider the family $\fr:=\{jz^5: j\in \N\}$ on $D$. Let $a(z)=e^z$, $\psi(z)=-e^z$. Then  for $n=0$ and $k=1$,  $f^{(k)}(z)-a(z)f^n(z)-\psi(z)$ has zeros of multiplicity $4$.  Clearly, all conditions of Theorem \ref{thm2} are satisfied but $\fr$ is not normal in $D$.
\end{example}
The next example shows that the condition  $\psi(z)\neq0$ can not be relaxed.
\begin{example}
Let $D:=\{z\in \C:|z|<1\}$. Consider the family $\fr:=\{jz^4: j\in \N\}$ on $D$. Let $a(z)=0$, $\psi(z)=z^{d-1}$. Then  for $k=1$,  $f^{(k)}(z)-a(z)f^n(z)-\psi(z)$ has zeros of multiplicity $3$.  Clearly, all conditions of Theorem \ref{thm2} are satisfied but $\fr$ is not normal in $D$.
\end{example}
The following example supports our Theorem \ref{thm2}.
\begin{example}
Let $D:=\{z\in \C:z\neq0\}$. Consider the family $\fr:=\{jz^4: j\in \N\}$ on $D$. Let $\displaystyle{a(z)=-{z^{-4}}}$, $\psi(z)=-z^{2}/4$. Then  for $n=2$ and $k=1$,  $f^{(k)}(z)-a(z)f^n(z)-\psi(z)$ has zeros of multiplicity $2$ in $D$ and all conditions of Theorem \ref{thm2} are satisfied. $\fr$ is normal in $D$.
\end{example}

We also investigate this situation for $f(z)+a(z)(f^{(k)})^n(z)-\psi(z)$ and prove the following result:
 \begin{theorem}\label{thm3}
Let $\fr$ be a family of meromorphic functions in a domain $D$. Let $d, k, m $ and $ n(\geq 2)$ be four positive integers such that ${\frac{nk+1}{(n-1)d}+ \frac{n(k+1)}{m(n-1)}< 1}$. Let $a(z)(\not\equiv 0) $ and $\psi(z)(\neq0)$ be  holomorphic functions in $D$. If for each $f\in\fr$, all zeros of $f$ are of multiplicity at least $d$  and all zeros of $f(z)+a(z)(f^{(k)})^n(z)-\psi(z)$ are of multiplicity at least $m$, then $\fr$ is normal in $D$.
\end{theorem}
For a family of holomorphic functions we have the following strengthened version:
\begin{theorem}\label{thm4}
Let $\fr$ be a family of holomorphic functions in a domain $D$. Let $d, k, m$ and $ n(\geq 2)$ be four positive integers such that ${\frac{nk+1}{nd}+ \frac{1}{m}< 1}$. Let $a(z)(\not\equiv 0)$ and $\psi(z)(\neq0)$ be  holomorphic functions in $D$. If for each $f\in\fr$, all zeros of $f$ are of multiplicity at least $d$  and all zeros of $f(z)+a(z)(f^{(k)})^n(z)-\psi(z)$ are of multiplicity at least $m$, then $\fr$ is normal in $D$.
\end{theorem}
The following example supports  Theorem \ref{thm4}.
\begin{example}
Let $D:=\{z\in \C:|z|>0\}$. Consider the family $\fr:=\{jz^4: j\in \N\}$ on $D$. Let $\displaystyle{a(z)={z^{-1}/16}}$, $\psi(z)=-z^{3}/4$. Then  for $k=1$ and $n=2$,  $f(z)+a(z)(f^{(k)})^n(z)-\psi(z)$ has zeros of multiplicity $2$ and all conditions of Theorem \ref{thm4} are satisfied. Clearly, $\fr$ is normal in $D$.
\end{example}

\section{Some Notation and results of Nevanlinna theory}
Let $\Delta=\{z: |z|<1\}$ be the unit disk. We use the following standard functions of value distribution theory, namely
\begin{center}
$T(r,f),  m(r,f),  N(r,f)\ \text{and}\ \overline{N}(r,f)$.
\end{center}
We let $S(r,f)$ be any function satisfying
\begin{center}
$S(r,f)=o\big(T(r,f)\big)$,  as $r\rightarrow +\ity,$
\end{center}
 possibly outside a set of finite measure.\\

 \begin{FFT}Let $f$ be a meromorphic function on $\C$ and $a$ be  a complex number. Then $$T\left(r, \frac{1}{f-a} \right)= T(r, f) + O(1). $$\end{FFT}

 \begin{LDT}Let $f$ be a non-constant meromorphic function on $\C$, and let $k$ be a positive integer. Then the equality $$m\left(r, \frac{f^{(k)}}{f}\right)=S(r, f)$$ holds for all $r\in [1, \ity)$ excluding a set of finite measure.\end{LDT}
 \section{Preliminary results}
        In order to prove our results we need  the following Lemmas. The well known Zalcman Lemma is a very important tool in the study of normal families. The following is a new version due to  Zalcman ~\cite{Zalc}.
\begin{lemma}\cite{Zalc, Zalc 1}\label{lem1}Let $\mathcal F$ be a family of meromorphic  functions in the unit disk  $\Delta$, with the property that for every function $f\in \mathcal F,$  the zeros of $f$ are of multiplicity at least $l$ and the poles of $f$ are of multiplicity at least $k$ . If $\mathcal F$ is not normal at $z_0$ in $\Delta$, then for  $-l< \alpha <k$, there exist
\begin{enumerate}
\item{ a sequence of complex numbers $z_n \rightarrow z_0$, $|z_n|<r<1$},
\item{ a sequence of functions $f_n\in \mathcal F$ },
\item{ a sequence of positive numbers $\rho_n \rightarrow 0$},
\end{enumerate}
such that $g_n(\zeta)=\rho_n^{\alpha}f_n(z_n+\rho_n\zeta) $ converges to a non-constant meromorphic function $g$ on $\C$ with $g^{\#}(\zeta)\leq g^{\#}(0)=1$. Moreover, $g$ is of order at most two. Here, $g^{\#}(z)=\frac{|g'(z)|}{1+|g(z)|^2}$ is the spherical derivative of $g$.
\end{lemma}

\begin{remark}\label{rem}
In Lemma \ref{lem1}, if $\fr$ is a family of holomorphic functions, then by Hurwitz's Theorem the limit function $g$ is a non-constant  entire function and the order of $g$ is at most 1.
\end{remark}
The following lemma is  Milloux's inequality.
\begin{lemma}\label{milloux}\cite{Hay, ccyang, Yang} \label{lemma 3}
Suppose $f(z)$ is a non-constant meromorphic function in the complex plane and $k$ is a positive integer. Then
\begin{equation}
T(r, f)< \overline{N}(r,f) + N\left(r, \frac{1}{f} \right) +\overline{ N}\left(r, \frac{1}{f^{(k)}-1} \right) -N\left(r, \frac{1}{f^{(k+1)}} \right) + S(r, f).
\end{equation}
\end{lemma}

Let $f$ be a non-constant meromorphic function in $\C$. A  differential polynomial $P$ of $f$ is defined by $\displaystyle {P(z):= \sum_{i=1}^{n}\al_{i}(z)\prod_{j=0}^{p}\left(f^{(j)}\left(z\right)\right)^{S_{ij}}},$ where $S_{ij}$'s are non-negative integers and    $\al_i(z)\not\equiv 0$ are small functions of $f$, that is $T(r,\al_i)=o\big(T(r,f)\big)$. The lower degree of the differential polynomial $P$ is defined by $$d(P):= \min_{1\leq i \leq n}\sum_{j=0}^{p}S_{ij} .$$

The following result was proved by Dethloff et al. in ~\cite{dethloff}.
\begin{lemma}\label{lem2}
Let $a_1, \ldots, a_q$ be distinct non-zero complex numbers.  Let $f$ be a non-constant  meromorphic function and let $P$  be a non-constant  differential polynomial of $f$ with $d(P)\geq 2.$ Then
\begin{equation}\notag
T(r,f)\leq \left(\frac{q\theta(P)+1}{qd(P)-1}\right)\overline{N}\left(r,\frac{1}{f}\right)+\frac{1}{qd(P)-1}\sum_{j=1}^{q}\overline{N}\left(r,\frac{1}{P-a_j}\right)+ S\left(r,f\right),
\end{equation}
for all $r\in [1,+\ity)$ excluding a set of finite Lebesgue measure, where $\displaystyle\theta(P):= \max_{1\leq i \leq n}\sum_{j=0}^{p}jS_{ij}.$\\

 Moreover, in the case of an entire function, we have
 \begin{equation}\notag
T(r,f)\leq \bigg(\frac{q\theta(P)+1}{qd(P)}\bigg)\overline{N}\left(r,\frac{1}{f}\right)+\frac{1}{qd(P)}\sum_{j=1}^{q}\overline{N}\left(r,\frac{1}{P-a_j}\right)+ S(r,f),
\end{equation}
for all $r\in [1,+\ity)$ excluding a set of finite Lebesgue measure.
\end{lemma}

This result was proved  by Hinchliffe in ~\cite{hin} for $q=1$.\\

\section{Proof of Main Results}

{\bf{Proof of Theorem \ref{thm1}}}: Since normality is a local property, we assume that
${D}=\Delta$. Suppose that $\fr$ is not normal
in $\Delta$. Then there exists at least one point $z_0$ such that
$\fr$ is not normal  at the point $z_0$ in $\Delta$. Without loss of
generality we assume that $z_0=0$. Then by Lemma \ref{lem1}, for $\al=-k$  there
exist
\begin{enumerate}
\item{ a sequence of complex numbers $z_j \rightarrow 0$, $|z_j|<r<1$},
\item{ a sequence of functions $f_j\in \mathcal F$},
\item{ a sequence of positive numbers $\rho_j \rightarrow 0$},
\end{enumerate}
such that $g_j(\zeta)=\rho_j^{-k}f_j(z_j+\rho_j\zeta) $ converges to a non-constant meromorphic function $g(\zeta)$ on $\C$. The zeros of $g(\zeta)$ are of multiplicity at least $d$ and has
poles with multiplicity at least $2$. Moreover, $g(\zeta)$ is of order at most 2. \\

We see that \begin{align*}
 f_j^{(k)}(z_j+\rho_j\zeta)&-a(z_j+\rho_j\zeta) f_j^n(z_j+\rho_j\zeta)-\psi(z_j+\rho_j\zeta)\\& = g_j^{(k)}(\zeta)- \rho_j^{nk} a(z_j+\rho_j\zeta)g_j({\zeta})-\psi(z_j+\rho_j\zeta)\\ & \rightarrow g^{(k)}(\zeta)-\psi(0).
\end{align*}

By Hurwitz's Theorem, we see that $g^{(k)}(\zeta)-\psi(0)$ has at least $m$ zeros. Now, by Milloux's inequality and Nevanlinna's First Fundamental Theorem we get
\begin{align*}
T(r,g)&\leq \overline{N}(r, g) + N\left(r, \frac{1}{g} \right) +\overline{N}\left(r, \frac{1}{g^{(k)}-\psi(0)} \right)-N\left(r, \frac{1}{g^{(k+1)}} \right) + S(r, g)\\
&\leq \frac{1}{2}N(r, g) + (k+1)\overline{N}\left(r, \frac{1}{g} \right)+ \overline{N}\left(r, \frac{1}{g^{(k)}-\psi(0)} \right)+ S(r, f) \\
&\leq  \frac{1}{2}N(r, g)+ \frac{k+1}{d}N\left(r, \frac{1}{g}\right) + \frac{1}{m}N\left(r, \frac{1}{g^{(k)}-\psi(0)} \right)+ S(r, g)\\
&\leq \left(\frac{1}{2}+\frac{k+1}{d} \right)T(r, g) +\frac{1}{m}\left(T(r, g) + k \overline{N}(r, g) \right) + S(r, g)\\
&\leq \left(\frac{1}{2}+\frac{k+1}{d}+\frac{1}{m}+\frac{k}{2m} \right)T(r, g) + S(r, g).
\end{align*}

Combining with assumption ${\frac{k+1}{d}+\frac{k+2}{2m} <\frac{1}{2}}$, we get $g$ is constant. This is a contradiction. Hence $\fr$ is a normal family, and this completes the proof of  Theorem \ref{thm1}. $\Box$\\

We can prove Theorem \ref{thm2} by the method of Theorem \ref{thm1}. Since $\fr$ is a family of holomorphic functions, so by Remark \ref{rem}, the limit function $g$  is a non-constant entire function with zeros of multiplicity at least $d$ and $g^{(k)}(\zeta)-\psi(0)$ has at least $m$ zeros. Now, by Lemma \ref{milloux} and Nevanlinna's First Fundamental Theorem we get
\begin{align*}
T(r,g)&\leq \overline{N}(r, g) + N\left(r, \frac{1}{g} \right) +\overline{N}\left(r, \frac{1}{g^{(k)}-\psi(0)} \right)-N\left(r, \frac{1}{g^{(k+1)}} \right) + S(r, g)\\
&\leq (k+1)\overline{N}\left(r, \frac{1}{g} \right)+ \overline{N}\left(r, \frac{1}{g^{(k)}-\psi(0)} \right)+ S(r, f) \\
&\leq   \frac{k+1}{d}N\left(r, \frac{1}{g}\right) + \frac{1}{m}N\left(r, \frac{1}{g^{(k)}-\psi(0)} \right)+ S(r, g)\\
&\leq \frac{k+1}{d} T(r, g) +\frac{1}{m}\left(T(r, g) + k \overline{N}(r, g) \right) + S(r, g)\\
&\leq \left(\frac{k+1}{d}+\frac{1}{m} \right)T(r, g) + S(r, g).
\end{align*}

Combining with assumption ${\frac{k+1}{d}+\frac{1}{m} < 1}$, we get $g$ is constant. This is a contradiction. Hence $\fr$ is a normal family, and this completes the proof of the Theorem \ref{thm2}. \hfill$\Box$\\

{\bf{Proof of Theorem \ref{thm3}}}:
Again we assume that
${D}=\Delta$. Suppose that $\fr$ is not normal
in $\Delta$. Then there exists at least one point $z_0$ such that
$\fr$ is not normal  at the point $z_0$ in $\Delta$. Without loss of
generality we assume that $z_0=0$. Then by Lemma \ref{lem1}, for $\al=-k$  there
exist
\begin{enumerate}
\item{ a sequence of complex numbers $z_j \rightarrow 0$, $|z_j|<r<1$},
\item{ a sequence of functions $f_j\in \mathcal F$},
\item{ a sequence of positive numbers $\rho_j \rightarrow 0$},
\end{enumerate}
such that $g_j(\zeta)=\rho_j^{-k}f_j(z_j+\rho_j\zeta) $ converges to a non-constant meromorphic function $g(\zeta)$ on $\C$. The zeros of $g(\zeta)$ are of multiplicity at least $d$ and has
poles with multiplicity at least $2$. Moreover, $g(\zeta)$ is of order at most 2. \\

We see that \begin{align*}
 f_j(z_j+\rho_j\zeta)&+a(z_j+\rho_j\zeta) \left(f_j^{(k)}\right)^n(z_j+\rho_j\zeta)-\psi(z_j+\rho_j\zeta)\\& = \rho_j^k g_j(\zeta)+ a(z_j+\rho_j\zeta) \left(g_j^{(k)}\right)^n({\zeta})-\psi(z_j+\rho_j\zeta)\\ & \rightarrow a(0) \left(g^{(k)}\right)^n(\zeta)-\psi(0).
\end{align*}

Let $a(0)=a$. By Hurwitz's Theorem, we see that $a \left(g^{(k)}\right)^n(\zeta)-\psi(0)$ has at least $m$ zeros. Now we invoke Lemma \ref{lem2} for the differential polynomial $P=a \left(g^{(k)}\right)^n(\zeta)$ and $q=1$. Note that, $d(P)=n (\geq 2)$ and $\theta(P)=nk.$ By Lemma \ref{lem2} and Nevanlinna's First Fundamental Theorem,
\begin{align*}
T(r,g)&\leq \frac{nk+1}{n-1}\overline{N}\left(r, \frac{1}{g} \right)+ \frac{1}{n-1}\overline{N}\left(r, \frac{1}{a\left(g^{(k)}\right)^n-\psi(0)} \right) +S(r, g)\\
&\leq\frac{nk+1}{(n-1)d}{N}\left(r, \frac{1}{g}\right)+ \frac{1}{(n-1)m}N\left(r, \frac{1}{a\left(g^{(k)}\right)^n-\psi(0)} \right) +S(r, g)\\
&\leq \frac{nk+1}{(n-1)d}{N}\left(r, \frac{1}{g} \right)+ \frac{1}{(n-1)m}T\left(r, \left(g^{(k)}\right)^n\right) +S(r, g)\\
&=\frac{nk+1}{(n-1)d}{N}\left(r, \frac{1}{g} \right)+ \frac{n}{(n-1)m}T\left(r, g^{(k)}\right) +S(r, g)\\
&=\frac{nk+1}{(n-1)d}{N}\left(r, \frac{1}{g} \right)+ \frac{n}{(n-1)m}\left(m\left(r, g^{(k)}\right)+N\left(r, g^{(k)}\right)\right) +S(r, g)\\
&\leq \frac{nk+1}{(n-1)d}{T}\left(r, g \right)+ \frac{n}{(n-1)m}\left(m\left(r, \frac{g^{(k)}}{g}\right)+m\left(r, g\right)\right)\\
& \qquad + \frac{n}{(n-1)m}\left(N(r, g) + k\overline{N}(r, g) \right) + S(r, g)\\
&\leq  \frac{nk+1}{(n-1)d}{T}\left(r, g \right)+ \frac{n}{(n-1)m}T(r, g)+ \frac{nk}{(n-1)m}T(r, g) +S(r, g)\\
&\leq \left(\frac{nk+1}{(n-1)d}+ \frac{n(1+k)}{(n-1)m} \right)T(r, g)+ S(r, g).
\end{align*}

Combining this with assumption $\displaystyle{{\frac{nk+1}{(n-1)d}+ \frac{n(1+k)}{(n-1)m} < 1}}$, we get that $g$ is a constant. This is a contradiction. Hence $\fr$ is a normal family and this completes the proof of  Theorem \ref{thm3}. \hfill$\Box$\\

Again, we can prove Theorem \ref{thm4} by the method of Theorem \ref{thm3}. Since $\fr$ is a family of holomorphic functions so by Remark \ref{rem}, the limit function $g$  is a non-constant entire function with zeros of multiplicity at least $d$ and $a \left(g^{(k)}\right)^n(\zeta)-\psi(0)$ has at least $m$ zeros. Now, we apply Lemma \ref{lem2} to the differential polynomial $P=a \left(g^{(k)}\right)^n(\zeta)$ and $q=1$. Note that, $d(P)=n (\geq 2)$ and $\theta(P)=nk.$ By Lemma \ref{lem2} and Nevanlinna's First Fundamental Theorem,
\begin{align*}
T(r,g)&\leq \frac{nk+1}{n}\overline{N}\left(r, \frac{1}{g} \right)+ \frac{1}{n}\overline{N}\left(r, \frac{1}{a\left(g^{(k)}\right)^n-\psi(0)} \right) +S(r, g)\\
&\leq\frac{nk+1}{nd}{N}\left(r, \frac{1}{g}\right)+ \frac{1}{nm}N\left(r, \frac{1}{a\left(g^{(k)}\right)^n-\psi(0)} \right) +S(r, g)\\
&\leq \frac{nk+1}{nd}{N}\left(r, \frac{1}{g} \right)+ \frac{1}{nm}T\left(r, \left(g^{(k)}\right)^n\right) +S(r, g)\\
&=\frac{nk+1}{nd}{N}\left(r, \frac{1}{g} \right)+ \frac{n}{nm}T\left(r, g^{(k)}\right) +S(r, g)\\
&=\frac{nk+1}{nd}{N}\left(r, \frac{1}{g} \right)+ \frac{1}{m}\left(m\left(r, g^{(k)}\right)+N\left(r, g^{(k)}\right)\right) +S(r, g)\\
&\leq \frac{nk+1}{nd}{T}\left(r, g \right)+ \frac{1}{m}\left(m\left(r, \frac{g^{(k)}}{g}\right)+m\left(r, g\right)\right)+ S(r, g)\\
&\leq  \frac{nk+1}{nd}{T}\left(r, g \right)+ \frac{1}{m}T(r, g)+ S(r, g)\\
&\leq \left(\frac{nk+1}{nd}+ \frac{1}{m} \right)T(r, g)+ S(r, g).
\end{align*}

Combining this with assumption $\displaystyle{\frac{nk+1}{nd}+ \frac{1}{m}< 1}$, we get $g$ is constant. This is a contradiction. Hence $\fr$ is a normal family. This completes the proof of the Theorem \ref{thm4}. \hfill$\Box$\\

{\bf{Acknowledgement:}}   The second author  is thankful to the faculty and the administrative unit of School  of  Mathematics,  Harish-Chandra  Research  Institute,  Allahabad  for  their warm hospitality during the preparation of this paper.

\end{document}